\documentclass[10pt,twosided]{article}
\usepackage[tbtags]{amsmath}
\usepackage{epsfig,amstext,amssymb,amsthm,latexsym}
\allowdisplaybreaks[4]
\usepackage{color}

\usepackage{amssymb,color}

\definecolor{c20}{rgb}{0.,0.7,0.}
\definecolor{c30}{rgb}{0.,0.,1.}
\definecolor{c40}{rgb}{1,0.1,0.7}
\definecolor{c50}{rgb}{1,0,0}
\definecolor{c60}{rgb}{1,0.9,0.1}

\def\eE#1{\textcolor{c30}{#1}}
\def\eE#1{#1}

\def\cE#1{\textcolor{c30}{#1}}
\def\cE#1{#1}
\def\aE#1{\textcolor{c30}{#1}}
\def\aE#1{#1}

\def\ccc#1{{\textcolor{c40}{#1}}}
\def\ccc#1{#1}
\def\wz#1{{\textcolor{c40}{#1}}}
\def\wz#1{#1}
\def\wzc#1{{\textcolor{c40}{#1}}}
\def\wzc#1{#1}

\usepackage{amssymb,color}

\newcommand{\kb}[1]{\boldsymbol{#1}}
\newcommand{\vk}[1]{\kb{#1}}

\newcommand{\abs}[1]{\lvert #1 \rvert}

\newcommand{\E}[1]{\mathbb{E}\left(#1\right)}

\newcommand{\pk}[1]{\mathbb{P} \left\{ #1 \right\} }

\newcommand{\R}{\!I\!\!R}

\newcommand{\N}{\!I\!\!N}
\newcommand{\inr}{\in \R}

\newcommand{\inn}{\in \N}

\newcommand{\ldot}{,\ldots,}

\newcommand{\limit}[1]{\lim_{#1 \to   \infty}}
\newcommand{\todis}{\stackrel{d}{\to}}

\newcommand{\BQN}{\begin{eqnarray}}
\newcommand{\EQN}{\end{eqnarray}}
\newcommand{\BQNY}{\begin{eqnarray*}}
\newcommand{\EQNY}{\end{eqnarray*}}

\newcommand{\BS}{\begin{sat}}
\newcommand{\ES}{\end{sat}}
\newcommand{\BT}{\begin{theo}}
\newcommand{\ET}{\end{theo}}
\newcommand{\BK}{\begin{korr}}
\newcommand{\EK}{\end{korr}}

\newcommand{\BD}{\begin{de}}
\newcommand{\ED}{\end{de}}

\newcommand{\BIT}{\begin{itemize}}
\newcommand{\EIT}{\end{itemize}}

\newcommand{\BRM}{\begin{remarks}}
\newcommand{\ERM}{\end{remarks}}

\newcommand{\BEL}{\begin{lem}}
\newcommand{\EEL}{\end{lem}}

\numberwithin{equation}{section}
\newtheorem{theo}{Theorem}[section]
\newtheorem{sat}[theo]{Proposition}
\newtheorem{de}[theo]{Definition}
\newtheorem{lem}{Lemma}[section]
\newtheorem{korr}[theo]{Corollary}

\newtheorem{remarks}[theo]{Remarks}

\newcommand{\prooftheo}[1]{ \textsc{Proof of Theorem} \ref{#1} }

\newcommand{\prooflem}[1]{\textsc{Proof of Lemma} \ref{#1}}

\newcommand{\QED}{\hfill $\Box$}

\newcommand{\EE}[1]{\mathbb{E}\left\{#1\right\}}

\def\E{\operatorname*{\mathbb{E}}}
\def\I{\operatorname*{\mathbb{I}}}

\newcommand{\netheo}[1]{{Theorem \ref{#1}}}

\def\bs{\boldsymbol}
\def\vps{\bs{\varepsilon}}

\newcommand{\nwc}{\newcommand}
\nwc{\COM}[1]{}

\def\IF{\infty}

\begin{document}
\centerline{\bf \Large \ccc{Joint \eE{Limiting} Distribution of} Minima and Maxima of Complete}

\centerline{\bf \Large and Incomplete Samples of Stationary Sequences}

\vspace{1cm}
\centerline{Enkelejd Hashorva and Zhichao Weng}

\centerline{University of Lausanne, Switzerland}


\vspace{0.81cm}
{\bf Abstract}: In \eE{the seminal contribution} \cite{Davis79} the joint weak convergence of maxima and minima of weakly dependent stationary sequences is derived under some mild asymptotic conditions. In this paper we address additionally the case of incomplete samples \aE{assuming that the average proportion
of incompleteness converges in probability to some random variable}. We show the joint weak convergence of the maxima and minima of both complete and incomplete samples. It turns out that for special cases, maxima and minima are asymptotically independent.

\eE{Berman condition}.


\section{Introduction}

Let $\{X_n, n\cE{\ge} 1\}$ be a stationary random sequence with marginal
distribution function $F$, i.e., all $X_i$'s have the same distribution function $F$, and $(X_{n+1} \ldot X_{n+j})$ has the same distribution as
$(X_{n+k+1} \ldot X_{n+k+j})$ for any $j,k,n \inn$. Suppose \aE{that} there exist sequences $a_n>0,\ c_n>0,\ b_n,\ d_n\inr$ and non-degenerate
distribution functions $G$ and $H$ such that (write next $\overline{F}=1- F, \overline{H}=1-H$)
\BQN\label{eqA}
\limit{n} F^n(u_n(x))=G(x), \qquad \limit{n} (\overline{F}(v_n(y)))^n= \overline{H}(y), 
\EQN
\COM{or equivalently
\begin{eqnarray}
\label{eq n F}
n\left(1-F\left(u_{n}(x)\right)\right)\to -\ln G(x),
\ nF\left(v_{n}(y)\right) \to -\ln \left(\overline{H}(y)\right)
\end{eqnarray}}
 where $u_{n}(x)=a_{n}x+b_{n}$ and $v_{n}(y)=c_{n}y+d_{n}$, $x,y \inr$.
In view of \cite{leadbetter1983extremes}
the normalized maxima $(\max_{i\le n} X_i - b_n)/a_n$ and the normalized minima $(\min_{i\le n} X_i - \aE{d_n)/c_n}$
are asymptotically independent if \eqref{eqA} holds and further \eE{$X_n,n\ge 1$} are independent.
For stationary random sequences \eE{the seminal contribution \cite{Davis79}  shows} that the joint limiting
behavior of the normalized \aE{maxima and minima}  is the same as that
\aE{in} the iid case, provided that some weak dependence conditions \aE{are imposed}.

The contributions \cite{kudrov2007maxima} and \cite{MR2307068}
studied  the asymptotic distribution of maxima of complete and incomplete samples. Several authors followed the aforementioned paper\eE{s} see e.g.,
\cite{MR2809873, MR2735150, MR2926583, MR2969062, HashorvaPengWeng} and the references therein.\\
We describe next the model which allows for the study of incomplete random samples
\aE{where some} of $X_i$'s can be observed whereas the others are not \eE{observable}.
Let $S_n=\sum_{i=1}^n \ccc{\varepsilon_i}$ denote the number of observed random variables from
$\{X_1,X_2,\cdots,X_n\}$, where $\varepsilon_i$ is the indicator of
the event that $X_i$ is observed. {Assume that for some constant} $\mathcal{P} \in [0,1]$ 
\BQN\label{eqeta}
\frac{S_n}{n}:=\frac{\sum_{i=1}^n \varepsilon_i}{n} \to \mathcal{P} , \text{
{in probability}} \EQN
{as $n\to \IF$}. \aE{Further,  \eE{set}}
\BQNY
M_n=\max_{1\le i\le n}X_{i}, \quad M_n(\bs{\varepsilon})=
\left\{
\begin{array}{ll}
\max\{X_{i}, \varepsilon_{i}=1, i\le n \}, & \mbox{if} \ \sum_{i=1}^n \varepsilon_{i}\ge 1;\\
\inf\{x|F(x)>0\}, & \mbox{otherwise}
\end{array}
\right.
\EQNY
and define similarly the minimum $m_n$ and $m_n(\bs{\varepsilon})$.\\
The seminal article \cite{MR2307068} {showed} that under some weak dependence
conditions (\eE{see for details Section 2 below})
\BQN \label{eq1.4}
\lim_{n\to\infty}\pk{M_{n}(\vps)\le a_n x+ b_n , M_{n}\le
a_n y+ b_n}=\mathcal{F}(x,y;\mathcal{P} )=: \ccc{G}^{\mathcal{P} }(x)G^{1-\mathcal{P} }(y) \EQN
holds for any   $x<y$. Recently \cite{MR2811020} showed that \eqref{eq1.4} still holds with
 $ \mathcal{F}(x,y;\mathcal{P} )=\EE{ G^{\mathcal{P} }(x)G^{1-\mathcal{P} }(y)}$, provided that \eqref{eqeta} \aE{is valid} with
 $\mathcal{P} $ being a random variable. A similar result addressing both maxima and minima when $\mathcal{P}$ is allowed to be random does
 not \eE{exist} in the literature. Therefore, in this paper we investigate the joint limiting distribution
of the normalized random vector
$(M_n(\bs{\vps}), m_n(\bs{\vps}), M_n, m_n)$ under some weak
dependence conditions similar to those given in \cite{Davis79} and \cite{MR2811020}
\aE{and assuming further that \eqref{eqeta} holds with $\mathcal{P}$ being random.}
Our result shows that $(M_n(\bs{\vps}),M_n)$ and $(m_n(\bs{\vps}), m_n)$ are asymptotically independent if the limit in \eqref{eqeta} holds with $
\mathcal{P} $ a \eE{deterministic} constant. Otherwise, if $\mathcal{P} $ is a non-degenerate random variable, this is not the case anymore. This fact is interesting \ccc{and} also expected since the incompleteness of the data influences both maxima and minima, and therefore the asymptotic independence is not always possible.

Brief organisation of the rest of the paper. In the next section we present our main result and then apply it to the interesting case of Gaussian stationary sequences. Proofs together with several auxiliary results are displayed in Section 3.

\section{Main Result}\label{sec2}
\cE{We shall consider below $\{X_n,n\ge 1\}$ a strictly stationary random sequence with marginal distribution function $F$.
In order to formulate our main result we need to present first the conditions of weak dependence for this random sequence which were given for
the case of complete random samples by Davis \cite{Davis79}.
Throughout in the sequel  $\{u_n(x_1), u_n(x_2), v_n(y_1), v_n(y_2), n\ge 1\}$ are given constants.}

{\bf Definition}: (Condition $D(u_n(x_1),
v_n(y_1), u_n(x_2), v_n(y_2))$ of \cite{Davis79}). We say that condition $D(u_n(x_1),
v_n(y_1), u_n(x_2), v_n(y_2))$ is satisfied, if for any $n$ and all
$A_1$, $A_2$, $B_1$, $B_2\subset \{1,2,\cdots,n\}$, such that $A_1\cap
A_2=\emptyset$, $B_1\cap B_2=\emptyset$ and $b-a\geq l_n$, where $a\in
A_1\cup A_2$ and $b\in B_1\cup B_2$ we have
\COM{\begin{eqnarray*}
&&
\left|\pk{\bigcap_{j\in A_1\cup B_1}\{X_j\leq
u_n(x_2)\}\cap\bigcap_{j\in A_2\cup B_2}\{X_j\leq
u_n(x_1)\}}\right.
\\
&&\qquad\quad
-\pk{\bigcap_{j \in A_1}\{X_j\leq
u_n(x_2)\}\cap\bigcap_{j \in
A_2}\{X_j\leq u_n(x_1)\}}
\times \left.\pk{\bigcap_{j \in B_1}\{X_j\leq
u_n(x_2)\}\cap\bigcap_{j \in B_2}\{X_j\leq
u_n(x_1)\}}\right|
\\
&
\leq&\alpha_{n,l_n};
\end{eqnarray*}
\begin{eqnarray*}
&&
\left|\pk{\bigcap_{j\in A_1\cup B_1}\{X_j>
v_n(y_2)\}\cap\bigcap_{j\in A_2\cup B_2}\{X_j>
v_n(y_1)\}}\right.
\\
&&\qquad
-\pk{\bigcap_{j \in A_1}\{X_j> v_n(y_2)\}\cap\bigcap_{j
\in A_2}\{X_j> v_n(y_1)\}}
\times \left. \pk{\bigcap_{j \in B_1}\{X_j>
v_n(y_2)\}\cap\bigcap_{j \in B_2}\{X_j> v_n(y_1)\}}\right|
\\
&
\leq&\alpha_{n,l_n};
\end{eqnarray*}
and
}
\begin{eqnarray*}
&&
\left|\pk{\bigcap_{j\in A_1\cup B_1}\{v_n(\widetilde{y}_2)<X_j\leq
u_n(\widetilde{x}_2)\}\cap\bigcap_{j\in A_2\cup B_2}\{v_n(\widetilde{y}_1)<X_j\leq
u_n(\widetilde{x}_1)\}}\right.
\\
&&\quad
-\pk{\bigcap_{j \in A_1}\{v_n(\widetilde{y}_2)<X_j\leq
u_n(\widetilde{x}_2)\}\cap\bigcap_{j \in A_2}\{v_n(\widetilde{y}_1)<X_j\leq
u_n(\widetilde{x}_1)\}}
\\
&&\quad\quad \times\left.\pk{\bigcap_{j \in
B_1}\{v_n(\widetilde{y}_2)<X_j\leq u_n(\widetilde{x}_2)\}\cap\bigcap_{j \in
B_2}\{v_n(\widetilde{y}_1)\le X_j\leq u_n(\widetilde{x}_1)\}}\right|
\leq \alpha_{n,l_n},
\end{eqnarray*}
where $\limit{n}\alpha_{n,l_n}=0$ for some sequence $l_n\to \infty$ with $l_n/n\to 0$ and
\cE{$\widetilde{x}_i = x_i+I\cdot \IF, \widetilde{y}_i =- (1- I)J\cdot \IF +y_i, i=1,2, I,J\in \{0,1\}$ (set $0 \cdot \IF:=0$).}

{\bf Definition}: (Condition $D^\prime(u_n(x),v_n(y))$ of \cite{Davis79}).
{\it 
For real sequences $\{u_n(x), v_n(y), n\ge 1\}$  condition $D^\prime(u_n(x),v_n(y))$ is satisfied if
\begin{eqnarray*}
&&
\limsup_{n\to \infty}n\sum^{[n/k]}_{j=1}\Big[\pk{X_1>u_{n}(x),X_{j+1}>u_{n}(x)}
+\pk{X_1>u_{n}(x),X_{j+1}\leq v_n(y)}
\\
&& \qquad\qquad
+ \pk{X_1\leq v_{n}(y),X_{j+1}>u_n(x)}
+\pk{X_1\leq v_n(y),X_{j+1}\leq v_{n}(y)}\Big] =o(1)
\end{eqnarray*}
as $k \to \infty$.
}

We state next our main result.
\BT
\label{th M m} Let $\{X_n, n\ge 1\}$ be a strictly stationary random
sequence with underlying distribution function $F$. Suppose that \eqref{eqA} holds for $u_{n}(x),v_{n}(y)$, $x,y \inr$.
Assume further that both condition  $D^\prime(u_n(x), v_n(y))$ \ccc{and} $D(u_n(x_1), v_n(y_1), u_n(x_2), v_n(y_2))$ hold for
$x_2<x_1$, $y_2>y_1$. If the indicator random
variables $\bs{\vps}=\{\varepsilon_n, n\ge 1\}$ are independent of $\{X_n,n\ge 1\}$ and \eqref{eqeta} is satisfies, then
\begin{eqnarray*}
&&
\lim_{n\to\infty} \pk{v_n(y_2)<m_n(\bs{\vps})\leq M_n(\bs{\vps})\leq
u_n(x_2),v_n(y_1)<m_n\leq M_n\leq u_n(x_1)}
\\
&&\qquad \qquad
=\EE{G^\mathcal{P}  (x_2)(\overline{H}(y_2))^\mathcal{P}  G^{1-\mathcal{P} }(x_1)(\overline{H}(y_1))^{1-\mathcal{P} }}.
\end{eqnarray*}
\ET


{\bf Remarks}: a) Under the conditions of Theorem \ref{th M m} for $y_1 <y_2$ we have
\begin{eqnarray*}
\lim_{n\to\infty}\pk{ m_n(\bs{\vps}) >v_n(y_2), m_n >v_n(y_1)}
=\EE{(\overline{H}(y_2))^\mathcal{P} (\overline{H}(y_1))^{1-\mathcal{P} }}.
\end{eqnarray*}
Further,
\begin{eqnarray*}
\lim_{n\to\infty}\pk{M_n(\bs{\vps})\leq
u_n(x_2), M_n\leq u_n(x_1)} =\EE{G^\mathcal{P} (x_2)G^{1-\mathcal{P} }(x_1)}
\end{eqnarray*}
holds with $x_2 < x_1$.\\
b) \wzc{Theorem \ref{th M m} implies}
\begin{eqnarray*}
\lim_{n\to\infty}\pk{v_n(y)<m_n(\bs{\vps})
\leq
M_n(\bs{\vps})\leq u_n(x)}=\EE{G^\mathcal{P} (x)(\overline{H}(y))^\mathcal{P} }.
\end{eqnarray*}
Hence, if $\mathcal{P} $ is a constant, then the maxima and minima are asymptotically independent.\\
c) Our result shows in particular the joint asymptotic convergence of $( m_n(\bs{\vps}), m_n)$ (and similarly for  $ (M_n(\bs{\vps}), M_n)$).
We have thus
$$ \eE{ \Biggl( \frac{m_n(\bs{\vps})- d_n}{c_n}, \frac{m_n- d_n}{c_n} \Biggr)} \todis (\mathcal{M}^*,\mathcal{M}), \quad n\to \IF$$
and consequently,
$$ \eE{ \Biggl( \frac{m_n(\bs{\vps})- m_n}{c_n}, \frac{m_n- d_n}{c_n} \Biggr)} \todis (\mathcal{M}^*-\mathcal{M} ,\mathcal{M}), \quad n\to \IF.$$
A similar result is given in \cite{kudrov2007maxima} when $\mathcal{P} $ is a \eE{deterministic} constant.

{\bf Example}. We consider next the case that $\{X_n,n\ge 1\}$ is a \eE{centered} stationary Gaussian sequence with correlation function $\rho_n=\EE{X_1X_{n+1}}<1,n\ge1$ such that $\EE{X_n}=0,\EE{X_n^2}=1,n\ge 1$. With the choice of constants
\BQN\label{eq1.2}
 a_n=c_n=1/\sqrt{2\ln n}, \quad b_n=d_n=\sqrt{2\ln
n}-\frac{\ln\ln n+ \ln 4\pi}{2 \sqrt{2\ln n}}
\EQN
condition \eqref{eqA} holds with \wzc{$u_n(x)=a_nx+b_n$, $v_n(y)=-c_ny-d_n$ and $\overline{H}(x)=\ccc{G}(x)=\exp(-\exp(-x))$} if further the Berman condition
\BQN\label{con:Be}
\limit{n} \rho_n \ln n =0
\EQN
is valid, see \cite{MR0161365, Berman92}. Note in passing that \eqref{eqA} also holds if
\BQN\label{con:Dav}
\sum_{n\ge 1} \abs{\ccc{\rho}_n}^p < \IF
\EQN for some $p>1$, see \cite{Davis79}. In view of \wzc{\cite{leadbetter1983extremes},
both condition $D(u_n(x_1),v_n(y_1),u_n(x_2),v_n(y_2))$ and $D^{\prime}(u_n(x),v_n(y))$ hold under \eqref{con:Be} or \eqref{con:Dav},
then the claim of \netheo{th M m} holds for such stationary Gaussian sequences}.

\section{Further Results and Proofs}\label{sec3}
\ccc{In order to prove the main theorem, we need some auxiliary results. Let
$\bs{\beta}=\{\cE{\beta_n, n\ge 1}\}$ be a non-random sequence taking values in $\{0,1\}$.}
Given some index set $I\subset \{1 \ldot  n\}$ we define
\begin{eqnarray*}
M(I,\bs{\beta})=
\left\{
\begin{array}{ll}
\max\{X_{i}, \beta_{i}=1, i\in I  \}, & \mbox{if} \ \sum_{i \in I}\beta_{i}\ge 1;\\
\inf\{x|F(x)>0\}, & \mbox{otherwise}
\end{array}
\right.
\end{eqnarray*}
and similarly for $m(I,\bs{\beta})$ where we consider instead of maximum, the minimum of $X_i$'s. If $J$ is another index set we
define $\tilde{d}(I,J):= \min_{i\in I, j\in J} \abs{i-j}$. Let $k$ be a fixed positive integer, $t=[n/k]$ and define
\BQNY
K_s=\{j: (s-1)t+1\le j\le st\}
\EQNY
for $1 \le s \le k$.
For \aE{a} random variable $\mathcal{P}\in [0,1]$ a.s.\ write
\begin{eqnarray*}
B_{r,k}=
\left\{\omega: \mathcal{P} (\omega) \in \Bigg\{
\begin{array}{ll}
[0,\frac{1}{2^k}], & r=0,\\
(\frac{r}{2^k},\frac{r+1}{2^k}], & 0< r \le 2^k-1\\
\end{array}
\right\}
\end{eqnarray*}
and \eE{set}
$$B_{r, k, \bs{\beta},n}=\{\omega: \varepsilon_{i}(\omega)=\beta_{i}, 1\le i\le n\}\cap B_{r, k}.$$

\BEL\label{le prod Is}
\ccc{\aE{If} condition $D(u_n(x_1),v_n(y_1),u_n(x_2),v_n(y_2))$ holds} \wzc{for $x_2<x_1$, $y_2>y_1$}, \aE{then} for $I_1,I_2,\ldots, I_k$ \eE{non-empty} subsets of $\{1,2,\cdots,n\}$  we have
\BQN\label{eqadd1}
&&\left|\pk{\bigcap_{s=1}^k v_n(y_2)<m(I_s,\bs{\beta})\le M(I_s, \bs{\beta})\le u_n(x_2), v_n(y_1)<m(I_s)\le M(I_s)\le u_n(x_1)}\right.\nonumber\\
&&\quad \left.-\prod_{s=1}^k \pk{v_n(y_2)<m(I_s,\bs{\beta})\le M(I_s, \bs{\beta})\le u_n(x_2), v_n(y_1)<m(I_s)\le M(I_s)\le u_n(x_1)}
\right|\nonumber\\
&\le&(k-1)\alpha_{n,l_n},
\EQN
provided that $ \min_{1 \le i < j \le k} \tilde{d}(I_i,I_j) \ge  l_n$.
\EEL

\prooflem{le prod Is}
For $k=2$, the inequality \eqref{eqadd1} is just the condition $D(u_n(x_1),v_n(y_1),u_n(x_2),v_n(y_2))$. Suppose
that inequality \eqref{eqadd1} holds for arbitrary $k-1$ set, such that the distance between any two of them is
not less then $l_n$. Define
$$
\mathcal{A}(I)=\{v_n(y_2)<m(I,\bs{\beta})\leq M(I,\bs{\beta})\leq u_n(x_2),v_n(y_1)<m(I)\leq M(I)\leq u_n(x_1)\}
$$
for any interval $I\in \{1,2,\cdots, n\}$. Using induction and the condition $D(u_n(x_1),v_n(y_1),u_n(x_2),v_n(y_2))$, we have
\begin{eqnarray*}
\left|P\left(\bigcap^k_{s=1}\mathcal{A}(I_s)\right) -
\prod^k_{s=1}P\left(\mathcal{A}(I_s)\right)\right|
&\le&\left|P\left(\bigcap^k_{s=1}\mathcal{A}(I_s)\right) -
P\left(\bigcap^{k-1}_{s=1}\mathcal{A}(I_s)\right)P\left(\mathcal{A}(I_k)\right)\right|\\
&&+\left|P\left(\bigcap^{k-1}_{s=1}\mathcal{A}(I_s)\right) -
\prod^{k-1}_{s=1}P\left(\mathcal{A}(I_s)\right)\right|P\left(\mathcal{A}(I_k)\right)\\
&\leq& \alpha_{n,l_n}+(k-2)\alpha_{n,l_n}=(k-1)\alpha_{n,l_n},
\end{eqnarray*}
thus the claim follows.
\QED

\BEL\label{le prod Ks}
Under the assumptions of Lemma \ref{le prod Is} we have
\BQNY
&&\big|\pk{v_n(y_2)<m_n(\bs{\beta})\le M_n(\bs{\beta}) \le u_n(x_2), v_n(y_1)<m_n\le M_n \le u_n(x_1)}\\
&&\quad -\prod_{s=1}^k\pk{v_n(y_2)<m(K_s,\bs{\beta})\le M(K_s,\bs{\beta}) \le u_n(x_2), v_n(y_1)<m(K_s)\le M(K_s) \le u_n(x_1)}\big|\\
&\le& (k-1)\alpha_{n,l_n}+(4k+2)l_n(\overline{F}(u_n(x_2))+F(v_n(y_2))).
\EQNY
\EEL

\prooflem{le prod Ks}
Using Lemma \ref{le prod Is}, the proof is similar to the proof of Lemma 3.3 in \cite{MR2809873}.
\QED

\def\psA{\Psi_n(x_1, y_1)}
\def\psB{\Psi_n(x_2, y_2)}
\def\FFF{ \aE{F_n(\vk{x},\vk{y})}}

\prooftheo{th M m}
Define in the following $\Psi_n(z_1,z_2)=\overline{F}(u_n(z_1))+F(v_n(z_2))$,
\BQNY
P(K_s,\bs{\varepsilon})=\pk{v_n(y_2)< m(K_s,\bs{\varepsilon})\le M(K_s,\bs{\varepsilon}) \le u_n(x_2),v_n(y_1)< m(K_s)\le M(K_s) \le u_n(x_1)}
\EQNY
\wzc{for $1\le s\le k$} and
\BQNY
P(n,\bs{\varepsilon})=\pk{v_n(y_2)< m_n(\bs{\varepsilon})\le M_n(\bs{\varepsilon}) \le u_n(x_2),v_n(y_1)< m_n\le M_n \le u_n(x_1)}.
\EQNY
Note that
\BQNY
&&\left|P(n,\bs{\varepsilon})-\EE{\prod_{s=1}^k\left[1-\frac{\mathcal{P}n\psB+(1-\mathcal{P})n\psA}{k}\right]}\right|\\
&\le&\sum^{2^k-1}_{r=0}\sum_{\bs{\beta}\in\{0,1\}^n}
\EE{\left|P(n,\bs{\beta})-\prod_{s=1}^k\left[1-\frac{\mathcal{P}n\psB+(1-\mathcal{P})n\psA}{k}\right]\right|\I(B_{r,k,\bs{\beta},n})}\\
&\le&E_1+E_2+E_3,
\EQNY
where
\BQNY
E_1=\sum^{2^k-1}_{r=0}\sum_{\bs{\beta}\in\{0,1\}^n}
\EE{\left|P(n,\bs{\beta})-\prod_{s=1}^k P(K_s,\bs{\beta})\right|\I(B_{r,k,\bs{\beta},n})},
\EQNY
\BQNY
E_2=\sum^{2^k-1}_{r=0}\sum_{\bs{\beta}\in\{0,1\}^n}
\EE{\left|\prod_{s=1}^k P(K_s,\bs{\beta})-
\prod_{s=1}^k\left[1-\frac{\frac{r}{2^k}n\psB+(1-\frac{r}{2^k})n\psA}{k}\right]\right|\I(B_{r,k,\bs{\beta},n})}
\EQNY
and
\BQNY
E_3&=&\sum^{2^k-1}_{r=0}\sum_{\bs{\beta}\in\{0,1\}^n}
\E\left\{\left|
\prod_{s=1}^k\left[1-\frac{\frac{r}{2^k}n\psB+(1-\frac{r}{2^k})n\psA}{k}\right]\right.\right.\\
&&\qquad\qquad\qquad\qquad\left.
\left.-\prod_{s=1}^k\left[1-\frac{\mathcal{P}n\psB+(1-\mathcal{P})n\psA}{k}\right]\right|\I(B_{r,k,\bs{\beta},n})\right\}.
\EQNY
\aE{Since
$$ \limit{n}n \psB= -\ln G(x_2)-\ln \overline{H}(y_2),$$
\wz{$\limit{n}\alpha_{n,l_n}=0$ and }$\limit{n} l_n/n=0$, then Lemma \ref{le prod Ks} implies}
\BQN\label{eqE1}
E_1\le(k-1)\alpha_{n,l_n}+(4k+2)\frac{l_n}{n}n\psB \to 0
\EQN
as $n \to \IF$. \aE{Next,  for $0\le r\le 2^k-1$}
\BQNY
&&\left[1-\frac{rt}{2^k}\psB+t\left(1-\frac{r}{2^k}\right)\psA\right]+\left[\frac{\sum_{j\in K_s}\beta_j}{t}-\frac{r}{2^k}\right]t \FFF\\
&\le&P(K_s,\bs{\beta})\\
&\le&\left[1-\frac{rt}{2^k}\psB+t\left(1-\frac{r}{2^k}\right)\psA\right]\\
&&+t\sum_{j=2}^t\pk{A_{s1},A_{sj}}+\left[\frac{\sum_{j\in K_s}\beta_j}{t}-\frac{r}{2^k}\right]t \FFF ,
\EQNY
where
$$ \FFF=F(u_n(x_2))-F(u_n(x_1))-  F(v_n(y_2))+F(v_n(y_1)) $$
and
$$A_{sj}=\{X_{(s-1)t+j}>u_n(x_2)\}\cup\{X_{(s-1)t+j}\le v_n(y_2)\}, \quad j=\{1,2,\ldots,t\}.$$
\aE{Hence, using Lemma 3 in \cite{MR2811020}}
\BQNY
E_2&\le&\sum^{2^k-1}_{r=0}\sum_{\bs{\beta}\in\{0,1\}^n}
 \sum_{s=1}^k \EE{\left|P(K_s, \bs{\beta})-\left[1-\frac{\frac{r}{2^k}n\psB+(1-\frac{r}{2^k})n\psA}{k}\right]\right|\I(B_{r,k,\bs{\beta},n})}\\
&\le&\sum^{2^k-1}_{r=0}\sum_{\bs{\beta}\in\{0,1\}^n}
 \sum_{s=1}^k \EE{\left|\sum_{j\in K_s}\frac{\beta_j}{t}-\frac{r}{2^k}\right|\frac{n}{k} \FFF \I(B_{r,k,\bs{\beta},n})}+n\sum_{j=2}^t\pk{A_{s1},A_{sj}}\\
&\le&\sum^{2^k-1}_{r=0}
\sum_{s=1}^k \EE{\left|\sum_{j\in K_s}\frac{\beta_j}{t}-\frac{r}{2^k}\right|\I(B_{r,k})}\frac{n}{k} \FFF
+n\sum_{j=2}^t\pk{A_{s1},A_{sj}}\\
&\le&\sum_{s=1}^k\left[2(2s-1)\left(d\left(\frac{S_{ts}}{ts},\mathcal{P} \right)+d\left(\frac{S_{t(s-1)}}{t(s-1)},\mathcal{P} \right)\right)+\frac{1}{2^k}\right]
\frac{n}{k} \FFF +n\sum_{j=2}^t\pk{A_{s1},A_{sj}},
\EQNY
where $d(X,Y)$ stands for Ky Fan metric, i.e., $d(X,Y)=\inf\{\epsilon, \pk{|X-Y|>\epsilon}<\epsilon\}$.
\aE{Since}  \wz{$\lim_{t\to \IF}d\left(\frac{S_{ts}}{ts},\mathcal{P} \right) = 0$} and
$$\lim_{n\to \IF}n \FFF =\ln G(x_2)-\ln G(x_1)+\ln \overline{H}(y_2)-\ln \overline{H}(y_1),$$
\wz{taking the limit as $n\to \IF$ and then as $t\to \IF$,}
we get
\BQN\label{eqE2}
\limsup_{n\to\IF}E_2 \le \frac{1}{2^k}[\ln G(x_2)-\ln G(x_1)+\ln \overline{H}(y_2)-\ln \overline{H}(y_1)]+ko(\frac{1}{k}).
\EQN
For $E_3$, we have
\BQN\label{eqE3}
E_3&\le&\sum^{2^k-1}_{r=0}\sum_{\bs{\beta}\in\{0,1\}^n}\sum_{s=1}^k
\EE{\left|\frac{r}{2^k}-\mathcal{P} \right|\frac{n}{k}(2-F(u_n(x_2))+F(v_n(y_2))-F(u_n(x_1))+F(v_n(y_1)))\I(B_{r,k,\bs{\beta},n})}\nonumber\\
&=&\sum^{2^k-1}_{r=0}
\EE{\left|\frac{r}{2^k}-\mathcal{P} \right|\I(B_{r,k})}n(2-F(u_n(x_2))+F(v_n(y_2))-F(u_n(x_1))+F(v_n(y_1)))\nonumber\\
&\le&\frac{1}{2^k}n(2-F(u_n(x_2))+F(v_n(y_2))-F(u_n(x_1))+F(v_n(y_1)))\nonumber\\
&\to&\frac{1}{2^k}[-\ln G(x_2)-\ln G(x_1)-\ln \overline{H}(y_2)-\ln \overline{H}(y_1)]
\EQN
as $n \to \IF$. Combining with \eqref{eqE1}-\eqref{eqE3}, we have
\BQNY
&&\left|\limsup_{n \to \IF}P(n,\bs{\varepsilon})-\E\left(1-\frac{-\ln G^{\mathcal{P} }(x_2)-\ln(\overline{H}(y_2))^{\mathcal{P} }-\ln G^{1-\mathcal{P} }(x_2)-\ln(\overline{H}(y_2))^{1-\mathcal{P} }}{k}\right)^k\right|\\
&\le& k o\left(\frac{1}{k}\right)+\frac{-\ln G(x_1)-\ln(\overline{H}(y_1))}{2^{k-1}}.
\EQNY
Hence, letting $k\to \IF$ yields the claim.
\QED

\COM{
For a fixed $s\in \{1,2,\ldots, k\}$ and the nonrandom sequence of $\bs{\beta}=\{\beta_n, n\ge 1\}\in \{0,1\}^N$,
let us consider the event
\BQNY
P(K_s,\bs{\beta})=\{v_n(y_2)< m(K_s,\bs{\beta})\le M(K_s,\bs{\beta}) \le u_n(x_2),v_n(y_1)< m(K_s)\le M(K_s) \le u_n(x_1)\}
\EQNY
and
\BQNY
P(n,\bs{\beta})=\{v_n(y_2)< m_n(\bs{\beta})\le M_n(\bs{\beta}) \le u_n(x_2),v_n(y_1)< m_n\le M_n \le u_n(x_1)\}.
\EQNY
For any $0\le r\le 2^k-1$ we have
\BQNY
&&\left[1-\frac{tr}{2^k}(\overline{F}(u_n(x_2))+F(v_n(y_2)))+t\left(1-\frac{r}{2^k}\right)(\overline{F}(u_n(x_1))+F(v_n(y_1)))\right]\\
&&+\left[\frac{\sum_{j\in K_s}\beta_j}{t}-\frac{r}{2^k}\right]t \FFF \\
&\le&P(K_s,\bs{\beta})\\
&\le&\left[1-\frac{tr}{2^k}(\overline{F}(u_n(x_2))+F(v_n(y_2)))+t\left(1-\frac{r}{2^k}\right)(\overline{F}(u_n(x_1))+F(v_n(y_1)))\right]\\
&&+t\sum_{j=2}^t\pk{A_{s1},A_{sj}}+\left[\frac{\sum_{j\in K_s}\beta_j}{t}-\frac{r}{2^k}\right]t \FFF ,
\EQNY
where $A_{sj}=\{X_{(s-1)t+j}>u_n(x_2)\}\cup\{X_{(s-1)t+j}\le v_n(y_2)\}$, $j=\{1,2,\ldots,t\}$.
Note that
\BQNY
&&\left|P(n, \bs{\beta})-\prod_{s=1}^k\left[1-\frac{\frac{r}{2^k}n(\overline{F}(u_n(x_2))+F(v_n(y_2)))+(1-\frac{r}{2^k})n(\overline{F}(u_n(x_1))+F(v_n(y_1)))}{k}\right]\right|\\
&\le&\left|P(n, \bs{\beta})-\prod_{s=1}^k P(K_s,\bs{\beta})\right|\\
&&+\left|\prod_{s=1}^k P(K_s,\bs{\beta})-\prod_{s=1}^k\left[1-\frac{\frac{r}{2^k}n(\overline{F}(u_n(x_2))+F(v_n(y_2)))+(1-\frac{r}{2^k})n(\overline{F}(u_n(x_1))+F(v_n(y_1)))}{k}\right]\right|\\
&=:&J_1+J_2.
\EQNY
Using Lemma \ref{le prod Ks}
\BQNY
J_1 \le (k-1)\alpha_{n,l_n}+(4k+2)(\overline{F}(u_n(x_2))+F(v_n(y_2))).
\EQNY
According to inequality
\BQN\label{eq prod to sum}
\left|\prod_{s=1}^k c_s- \prod_{s=1}^k e_s \right| \le \sum_{s=1}^k|c_s-e_s|
\EQN
holds for all $c_s,e_s\in [0,1]$, and since $0\le 1-\frac{tr}{2^k}(\overline{F}(u_n(x_2))+F(v_n(y_2)))+t(1-\frac{r}{2^k})(\overline{F}(u_n(x_1))+F(v_n(y_1))) \le 1$, we have
\BQNY
J_2&\le& \sum_{s=1}^k\left|P(K_s,\bs{\beta})-\left[1-\frac{\frac{r}{2^k}n(\overline{F}(u_n(x_2))+F(v_n(y_2)))+(1-\frac{r}{2^k})n(\overline{F}(u_n(x_1))+F(v_n(y_1)))}{k}\right]\right|\\
&\le& n \sum_{j=2}^t\pk{A_{11},A_{1j}}+\sum_{s=1}^k\frac{\left|\sum_{j\in K_s}\frac{\beta_j}{t}-\frac{r}{2^k}\right|}{k}n \FFF .
\EQNY
Furthermore, using \eqref{eq prod to sum} again,
\BQNY
&&\E \sum_{\bs{\beta}\in \{0,1\}^n}\left|\prod_{s=1}^k\left[1-\frac{\frac{r}{2^k}n(\overline{F}(u_n(x_2))+F(v_n(y_2)))+(1-\frac{r}{2^k})n(\overline{F}(u_n(x_1))+F(v_n(y_1)))}{k}\right]\right.\\
&&\quad\left.-\prod_{s=1}^k\left[1-\frac{\mathcal{P}  n(\overline{F}(u_n(x_2))+F(v_n(y_2)))+(1-\mathcal{P} )n(\overline{F}(u_n(x_1))+F(v_n(y_1)))}{k}\right]\right| \I(B_{r,k,\bs{\beta}})\\
&\le&\sum_{s=1}^k\E \left|\frac{r}{2^k}-\mathcal{P} \right|\I(B_{r,k})\frac{n(2-F(u_n(x_2))+F(v_n(y_2))-F(u_n(x_1))+F(v_n(y_1)))}{k}\\
&\le&\frac{n(2-F(u_n(x_2))+F(v_n(y_2))-F(u_n(x_1))+F(v_n(y_1)))}{2^k}\pk{B_{r,k}}.
\EQNY
Since $\{X_n\}$, $\{\bs{\vps}_n\}$ and $\mathcal{P} $ are independence, we have
\BQNY
\sum_{\bs{\beta}\in \{0,1\}^n}\E P(n, \bs{\beta})\I(B_{r,k,\bs{\beta}})=P(n,\bs{\vps})\pk{B_{r,k}}.
\EQNY
Hence we have
\BQN\label{eq J r k}
J_{r,k}&=&\sum_{\bs{\beta}\in \{0,1\}^n}\E\left|P(n,\bs{\vps})
-\prod_{s=1}^k\left[1-\frac{\mathcal{P}  n(\overline{F}(u_n(x_2))+F(v_n(y_2)))+(1-\mathcal{P} )n(\overline{F}(u_n(x_1))+F(v_n(y_1)))}{k}\right]\right|\I(B_{r,k,\bs{\beta}})\nonumber\\
&\le&((k-1)\alpha_{n,l_n}+(4k+2)(\overline{F}(u_n(x_2))+F(v_n(y_2))))\pk{B_{r,k}}+ n \sum_{j=2}^t\pk{A_{11},A_{1j},B_{r,k}}\nonumber\\
&&+\E \sum_{s=1}^k\frac{\left|\sum_{j\in K_s}\frac{\vps_j}{t}-\frac{r}{2^k}\right|}{k}n \FFF \I(B_{r,k})\nonumber\\
&&+\frac{n(2-F(u_n(x_2))+F(v_n(y_2))-F(u_n(x_1))+F(v_n(y_1)))}{2^k}\pk{B_{r,k}}.
\EQN
Using Lemma 3 in Krajka (2011) we have
\BQNY
&&\sum_{r=0}^{2^k-1}\E \left|\sum_{j\in K_s}\frac{\vps_j}{t}-\frac{r}{2^k}\right|\I(B_{r,k}) \\
&\le&\E \left|\sum_{j\in K_s}\frac{\vps_j}{t}-\mathcal{P} \right|+\sum_{r=0}^{2^k-1}\E\left|\mathcal{P} -\frac{r}{2^k}\right|\I(B_{r,k})\\
&\le&\E \left|\frac{S_{ts}-S_{t(s-1)}}{t}-\mathcal{P} \right|+\frac{1}{2^k}\\
&\le& 2(2s-1)\left(d\left(\frac{S_{ts}}{ts},\mathcal{P} \right)+d\left(\frac{S_{t(s-1)}}{t(s-1)},\mathcal{P} \right)\right)+\frac{1}{2^k},
\EQNY
where $d(X,Y)$ stands for Ky Fan metric i.e., $d(X,Y)=\inf\{\epsilon, \pk{|X-Y|>\epsilon}<\epsilon\}$.
Thus taking a sum $\sum_{r=0}^{2^k-1}$ of the left- and right-hand side of \eqref{eq J r k} we have
\BQNY
\sum_{r=0}^{2^k-1}J_{r,k} &\le& (k-1)\alpha_{n,l_n}+(4k+2)(\overline{F}(u_n(x_2))+F(v_n(y_2)))+n \sum_{j=2}^t\pk{A_{11},A_{1j}}\\
&&+\sum_{s=1}^k\left[2(2s-1)\left(d\left(\frac{S_{ts}}{ts},\mathcal{P} \right)+d\left(\frac{S_{t(s-1)}}{t(s-1)},\mathcal{P} \right)\right)+\frac{1}{2^k}\right]
\frac{n \FFF }{k}\\
&&+\frac{n(2-F(u_n(x_2))+F(v_n(y_2))-F(u_n(x_1))+F(v_n(y_1)))}{2^k}.
\EQNY
Using \eqref{eqA}, condition $D^\prime(u_n(x_2),v_n(y_2))$ and $\lim_{t\to \IF}d\left(\frac{S_{ts}}{ts},\mathcal{P} \right)=0$, we have
\BQNY
&&\left|\lim_{n \to \IF}P(n,\bs{\vps})-\E\left(1-\frac{-\ln G^{\mathcal{P} }(x_2)-\ln(\overline{H}(y_2))^{\mathcal{P} }-\ln G^{1-\mathcal{P} }(x_2)-\ln(\overline{H}(y_2))^{1-\mathcal{P} }}{k}\right)^k\right|\\
&\le& k o\left(\frac{1}{k}\right)+\frac{-\ln G(x_1)-\ln(\overline{H}(y_1))}{2^{k-1}}.
\EQNY
Hence, letting $k\to \IF$ yields the claim. \QED}

{\textbf{Acknowledgments.}}
E.\ Hashorva acknowledges support from
the Swiss National Science Foundation grant 200021-1401633/1; Z. Weng has been partially
supported by the Swiss National Science Foundation Project under
grant 200021-134785
and by the project RARE -318984 (a Marie Curie
IRSES Fellowship within the 7th European Community Framework Programme).

\bibliographystyle{plain}
\bibliography{gausMinMaxM}

\newcommand{\nosort}[1]{}\def\polhk#1{\setbox0=\hbox{#1}{\ooalign{\hidewidth
  \lower1.5ex\hbox{`}\hidewidth\crcr\unhbox0}}}
  \def\polhk#1{\setbox0=\hbox{#1}{\ooalign{\hidewidth
  \lower1.5ex\hbox{`}\hidewidth\crcr\unhbox0}}} \def\cprime{$'$}
  \def\cprime{$'$} \def\cprime{$'$}
\begin{thebibliography}{10}

\bibitem{MR0161365}
S.M. Berman.
\newblock Limit theorems for the maximum term in stationary sequences.
\newblock {\em Ann. Math. Statist.}, 35:502--516, 1964.

\bibitem{Berman92}
S.M. Berman.
\newblock {\em Sojourns and extremes of stochastic processes}.
\newblock The Wadsworth \& Brooks/Cole Statistics/Probability Series. Wadsworth
  \& Brooks/Cole Advanced Books \& Software, Pacific Grove, CA, 1992.

\bibitem{MR2735150}
L.~Cao and Z.~Peng.
\newblock Asymptotic distributions of maxima of complete and incomplete samples
  from strongly dependent stationary {G}aussian sequences.
\newblock {\em Appl. Math. Lett.}, 24(2):243--247, 2011.

\bibitem{Davis79}
R.A. Davis.
\newblock Maxima and minima of stationary sequences.
\newblock {\em Ann. Probab.}, 7(3):453--460, 1979.

\bibitem{HashorvaPengWeng}
E.~Hashorva, Z.~Peng, and Z.~Weng.
\newblock On {P}iterbarg theorem for maxima of stationary {G}aussian sequences.
\newblock {\em Lithuanian Mathematical Journal}, in press, 2013.

\bibitem{MR2811020}
T.~Krajka.
\newblock The asymptotic behaviour of maxima of complete and incomplete samples
  from stationary sequences.
\newblock {\em Stochastic Process. Appl.}, 121(8):1705--1719, 2011.

\bibitem{kudrov2007maxima}
A.V. Kudrov and V.I. Piterbarg.
\newblock On maxima of partial samples in {G}aussian sequences with
  pseudo-stationary trends.
\newblock {\em Lithuanian Mathematical Journal}, 47(1):48--56, 2007.

\bibitem{leadbetter1983extremes}
M.R. Leadbetter, G.~Lindgren, and H.~Rootz{\'e}n.
\newblock {\em Extremes and related properties of random sequences and
  processes}, volume~11.
\newblock Springer Verlag, 1983.

\bibitem{MR2307068}
P.~Mladenovi{\'c} and V.I. Piterbarg.
\newblock On asymptotic distribution of maxima of complete and incomplete
  samples from stationary sequences.
\newblock {\em Stochastic Process. Appl.}, 116(12):1977--1991, 2006.

\bibitem{MR2926583}
Z.~Peng, B.~Tong, and S.~Nadarajah.
\newblock Almost sure central limit theorems of the partial sums and maxima
  from complete and incomplete samples of stationary sequences.
\newblock {\em Stoch. Dyn.}, 12(3):1150026, 19, 2012.

\bibitem{MR2809873}
Z.~Peng, Z.~Weng, and S.~Nadarajah.
\newblock Joint limiting distributions of maxima and minima for complete and
  incomplete samples from weakly dependent stationary sequences.
\newblock {\em J. Comput. Anal. Appl.}, 13(5):875--880, 2011.

\bibitem{MR2969062}
Z.~Tan and Y.~Wang.
\newblock Some asymptotic results on extremes of incomplete samples.
\newblock {\em Extremes}, 15(3):319--332, 2012.

\end{thebibliography}

\end{document}